\documentclass[11pt,fleqn]{article}
\usepackage{amssymb}
\usepackage{amsmath}
\usepackage{graphicx}

\begin{document}
\title{Analytical computation of boundary integrals for the Helmholtz equation in three dimensions}
\author{Nail A. Gumerov and Ramani Duraiswami \\
Perceptual Interfaces and Reality Lab, UMIACS,\\
 University of Maryland, College Park, MD} 
\maketitle
\begin{abstract}
A key issue in the solution of partial differential equations via integral equation methods is the evaluation of possibly singular integrals involving the Green's function and its derivatives multiplied by simple functions over discretized representations of the boundary. For the Helmholtz equation, while many authors use numerical quadrature to evaluate these boundary integrals, we present analytical expressions for such integrals over flat polygons in the form of infinite series. These can be efficiently truncated based on the accurate error bounds, which is key to their integration in methods such as the Fast Multipole Method.
\end{abstract}

\section{Introduction}
Boundary integral methods can be used to solve the Helmholtz equation either in the direct form, e.g., \cite{Gumerov2009:JASA}, or the indirect form e.g., \cite{Gumerov2013:POMA}.
The boundary in these approaches was discretized with flat triangles or quadrilaterals, and the discrete integrals evaluated. In this brief note, we provide expressions for computation of boundary integrals for the Helmholtz equation in three dimensions. Such integrals can be used both in the direct and indirect boundary element methods (BEM). Particularly, in our paper \cite{Gumerov2013:POMA} we did not provide the expressions for the integrals due to the lack of space for short communications

Given surface $S$, the problem is to compute the single and double layer potentials and their derivatives,
\begin{eqnarray}
L_{k}\left[ \sigma \right] \left( \mathbf{r}\right)  &=&\int_{S}\sigma
\left( \mathbf{r}^{\prime }\right) G_{k}\left( \mathbf{r},\mathbf{r}^{\prime
}\right) dS\left( \mathbf{r}^{\prime }\right) ,  \label{int1} \\
M_{k}\left[ \mu \right] \left( \mathbf{r}\right)  &=&\int_{S}\mu \left( 
\mathbf{r}^{\prime }\right) \frac{\partial G_{k}\left( \mathbf{r},\mathbf{r}%
^{\prime }\right) }{\partial n\left( \mathbf{r}^{\prime }\right) }dS\left( 
\mathbf{r}^{\prime }\right) ,  \nonumber \\
\mathbf{L}_{k}^{\prime }\left[ \sigma \right] \left( \mathbf{r}\right) 
&=&\nabla \int_{S}\sigma \left( \mathbf{r}^{\prime }\right) G_{k}\left( 
\mathbf{r},\mathbf{r}^{\prime }\right) dS\left( \mathbf{r}^{\prime }\right) ,
\nonumber \\
\mathbf{M}_{k}^{\prime }\left[ \mu \right] \left( \mathbf{r}\right) 
&=&\nabla \int_{S}\mu \left( \mathbf{r}^{\prime }\right) \frac{\partial
G_{k}\left( \mathbf{r},\mathbf{r}^{\prime }\right) }{\partial n\left( 
\mathbf{r}^{\prime }\right) }dS\left( \mathbf{r}^{\prime }\right) . 
\nonumber
\end{eqnarray}%
where $\sigma $ and $\mu $ are the single and double layer densities, $\mathbf{n}$ is the normal to the surface, and $G_{k}\left( \mathbf{r}, \mathbf{r}^{\prime }\right) $ is the free-space Green function for given wavenumber $k$,     
\begin{equation}
G_{k}\left( \mathbf{r},\mathbf{r}^{\prime }\right) =\frac{e^{ik\left|  \mathbf{r}-\mathbf{r}^{\prime }\right| }}{4\pi \left| \mathbf{r}-\mathbf{r}^{\prime }\right| }.  \label{int2}
\end{equation}
\begin{figure}[tbh]
	\begin{center}
		\includegraphics[width=0.96\textwidth, trim=0.15in 1.75in 0.5in	0.1in]{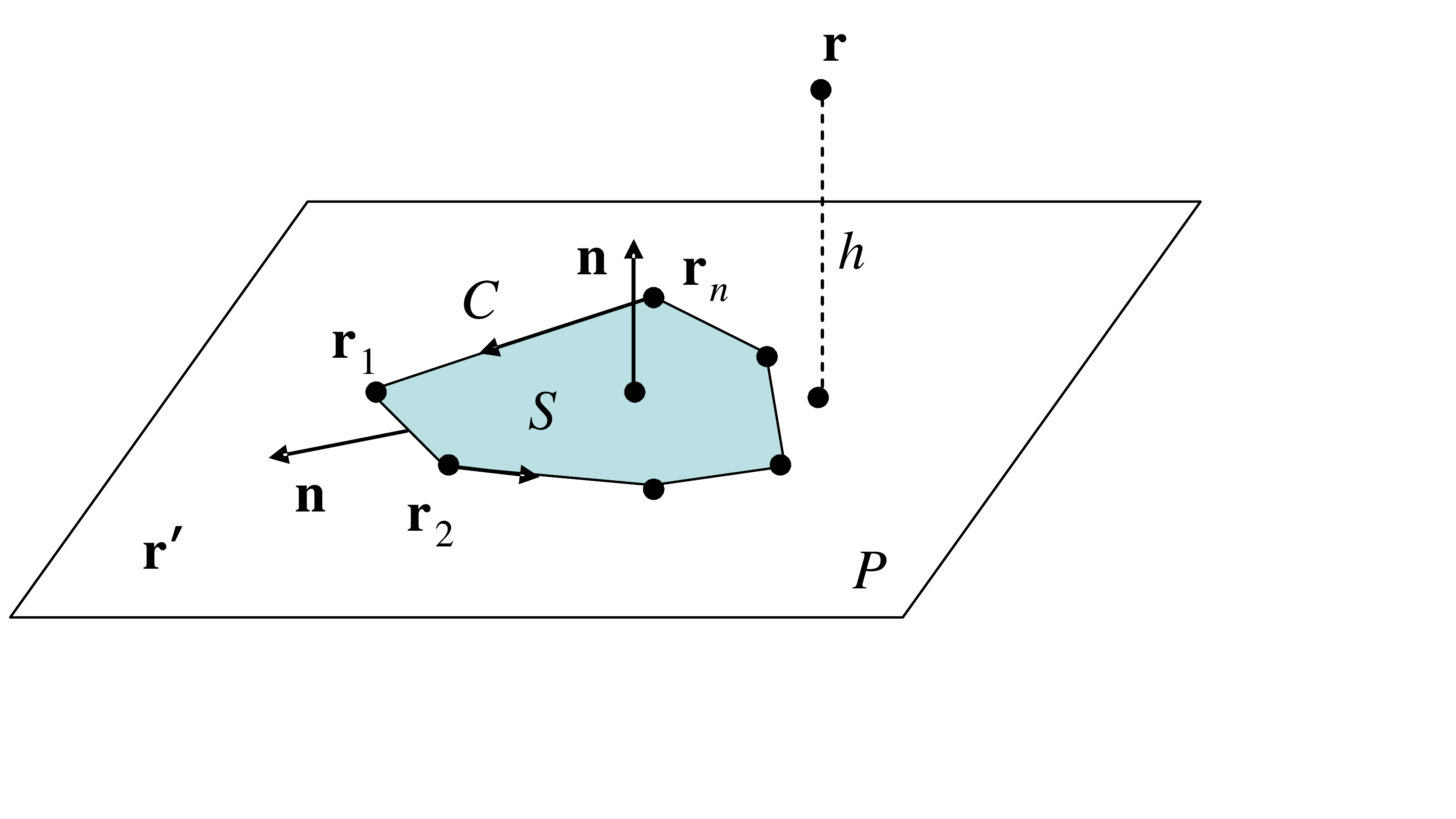}
	\end{center}
	\caption{Boundary integrals over flat polygons: notation. The polygon is defined by an ordered tuple of points $\left({\bf r}_1, {\bf r}_2, \ldots, {\bf r}_n, {\bf r}_1\right)$, which are arranged according to the right hand rule so that the normal $\textbf{n}$ points outward. The approach to computing the integrals uses the Gauss' divergence theorem and reduces the surface integral to a contour integral which is evaluated over the line segments. }
	\label{Fig1}
\end{figure}

\subsection{Computation of boundary integrals}
In center panel collocation with the boundary discretized by triangles, the integrals are approximated by $N$ integrals over flat patches (e.g.,
triangles $\Delta _{l}$) centered at $\mathbf{r}_{l}^{(c)}$, 
\begin{eqnarray}
L_{k}\left[ \sigma \right] \left( \mathbf{r}\right)  &=&\sum_{l^{\prime
}=1}^{N}\int_{\Delta _{l^{\prime }}}\sigma \left( \mathbf{r}^{\prime
}\right) G_{k}\left( \mathbf{r},\mathbf{r}^{\prime }\right) dS\left( \mathbf{%
r}^{\prime }\right) \approx \sum_{l^{\prime }=1}^{N}L_{ll^{\prime }}\sigma
_{l^{\prime }},  \label{i1} \\
\sigma _{l^{\prime }} &=&\sigma \left( \mathbf{r}_{l^{\prime }}^{(c)}\right)
,\quad L_{ll^{\prime }}=\int_{\Delta _{l^{\prime }}}G_{k}\left( \mathbf{r},%
\mathbf{r}^{\prime }\right) dS\left( \mathbf{r}^{\prime }\right) .  \nonumber
\end{eqnarray}%
Similar expressions can be written for other boundary integrals. For a flat
patch $\Delta $ with vertices ($\mathbf{r}_{1},\mathbf{...},\mathbf{r}_{n})$
defining the polygon as a right-hand oriented contour $C$ with normal $%
\mathbf{n}$ (see Fig. \ref{Fig1}), we need to compute the integrals 
\begin{eqnarray}
L_{k\Delta }\left( \mathbf{r}\right)  &=&\int_{\Delta }G_{k}\left( \mathbf{%
r,r}^{\prime }\right) dS\left( \mathbf{r}^{\prime }\right) ,\quad M_{k\Delta
}\left( \mathbf{r}\right) =\int_{\Delta }\mathbf{n}\cdot \nabla _{\mathbf{r}%
^{\prime }}G_{k}\left( \mathbf{r}^{\prime }\mathbf{,r}\right) dS\left( 
\mathbf{r}^{\prime }\right) ,  \label{1} \\
\mathbf{L}_{k\Delta }^{\prime }\left( \mathbf{r}\right)  &=&\nabla _{\mathbf{%
r}}\int_{\Delta }G_{k}\left( \mathbf{r}^{\prime }\mathbf{,r}\right) dS\left( 
\mathbf{r}^{\prime }\right) ,\quad \mathbf{M}_{k\Delta }^{\prime }\left( 
\mathbf{r}\right) =\nabla _{\mathbf{r}}\int_{\Delta }\mathbf{n}\cdot \nabla
_{\mathbf{r}^{\prime }}G_{k}\left( \mathbf{r}^{\prime }\mathbf{,r}\right)
dS\left( \mathbf{r}^{\prime }\right) .  \nonumber
\end{eqnarray}%
Note that%
\begin{gather}
G_{k}\left( \mathbf{r,r}^{\prime }\right) =\nabla _{s}^{\prime }\cdot \left[ 
\mathbf{\rho }f_{k}\left( \rho ;h\right) \right] ,\quad f_{k}\left( \rho
;h\right) =\frac{e^{ikr}-e^{ik\left| h\right| }}{4\pi ik\rho ^{2}},
\label{2} \\
h=\mathbf{n\cdot r,\quad \rho =r}^{\prime }-\mathbf{r+n}h,\quad r=\left| 
\mathbf{r}-\mathbf{r}^{\prime }\right| =\sqrt{\rho ^{2}+h^{2}},\quad \rho
=\left| \mathbf{\rho }\right| ,  \nonumber
\end{gather}%
where $\nabla _{s}^{\prime }$ is the surface gradient operator with respect
to $\mathbf{r}^{\prime }$ and $h$ is the projection of $\mathbf{r}$ to the
normal. Using the Gauss divergence theorem, we obtain%
\begin{equation}
L_{k\Delta }\left( \mathbf{r}\right) =\sum_{j=1}^{n}I_{j},\quad
I_{j}=\int_{C_{j}}\left( \mathbf{n}_{j}^{\prime }\cdot \mathbf{\rho }\right)
f_{k}\left( \rho ;h\right) dC,  \label{3}
\end{equation}%
where $C_{j}$ are line segments constituting the triangle sides and $\mathbf{%
n}_{j}^{\prime }\mathbf{\ }$are the normals to these segments outward to
the triangle. Expressions for $I_{j}$ can be simplified to%
\begin{eqnarray}
I_{j} &=&H\left( l_{j}-x^{\prime },y^{\prime },z^{\prime }\right) -H\left(
-x^{\prime },y^{\prime },z^{\prime }\right) ,\quad j=1,2,3,  \label{4} \\
H\left( x,y^{\prime },z^{\prime }\right)  &=&-z^{\prime }\int f\left( \rho
;y^{\prime }\right) dx=\frac{-z^{\prime }}{4\pi ik}\int \frac{%
e^{ikr}-e^{iky^{\prime }}}{x^{2}+z^{\prime 2}}dx,\quad r=\sqrt{%
x^{2}+y^{\prime 2}+z^{\prime 2}},\quad y^{\prime }=\left| h\right| , 
\nonumber
\end{eqnarray}%
where $l_{j}$ is the length of segment $C_{j}$, $\left( x^{\prime
},y^{\prime },z^{\prime }\right) $ are the coordinates of the evaluation
point $\mathbf{r}$ in the local right-hand oriented Euclidean reference
frame centered at the triangle vertex, which $x$-axis (basis vector $\mathbf{%
i}_{x}^{\prime }$) is directed along the segment, and the other basis
vectors are $\mathbf{i}_{y}^{\prime }=\mathbf{n}$, $\mathbf{i}_{z}^{\prime }=%
\mathbf{n}_{j}^{\prime }$. Expressions for the surface integrals can be
obtained by summation of the line integrals: 
\begin{eqnarray}
L_{k\Delta }\left( \mathbf{r}\right)  &=&\sum_{j=1}^{n}I_{j}\left(
x_{j}^{\prime },h,z_{j}^{\prime }\right) ,\quad I_{j}\left( x^{\prime
},y^{\prime },z^{\prime }\right) =H\left( l_{j}-x^{\prime },y^{\prime
},z^{\prime }\right) -H\left( -x^{\prime },y^{\prime },z^{\prime }\right) ,
\label{4.1} \\
\mathbf{L}_{k\Delta }^{\prime }\left( \mathbf{r}\right)  &=&\sum_{j=1}^{n}%
\mathbf{I}_{j}^{\prime }\left( x_{j}^{\prime },h,z_{j}^{\prime }\right)
,\quad \mathbf{I}_{j}^{\prime }\left( x^{\prime },y^{\prime },z^{\prime
}\right) =\mathbf{H}^{\prime }\left( l_{j}-x^{\prime },y^{\prime },z^{\prime
}\right) -\mathbf{H}^{\prime }\left( -x^{\prime },y^{\prime },z^{\prime
}\right) ,\quad   \nonumber \\
M_{k\Delta }\left( \mathbf{r}\right)  &=&\sum_{j=1}^{n}J_{j}\left(
x_{j}^{\prime },h,z_{j}^{\prime }\right) ,\quad J_{j}\left( x^{\prime
},y^{\prime },z^{\prime }\right) =K\left( l_{j}-x^{\prime },y^{\prime
},z^{\prime }\right) -K\left( -x^{\prime },y^{\prime },z^{\prime }\right) , 
\nonumber \\
\mathbf{M}_{k\Delta }^{\prime }\left( \mathbf{r}\right)  &=&\sum_{j=1}^{n}%
\mathbf{J}_{j}^{\prime }\left( x_{j}^{\prime },h,z_{j}^{\prime }\right)
,\quad \mathbf{J}_{j}^{\prime }\left( x^{\prime },y^{\prime },z^{\prime
}\right) =\mathbf{K}^{\prime }\left( l_{j}-x^{\prime },y^{\prime },z^{\prime
}\right) -\mathbf{K}^{\prime }\left( -x^{\prime },y^{\prime },z^{\prime
}\right) ,  \nonumber
\end{eqnarray}%
where 
\begin{equation}
\mathbf{H}^{\prime }=-\mathbf{i}_{x}^{\prime }\frac{\partial H}{\partial x}+%
\mathbf{i}_{y}^{\prime }\frac{\partial H}{\partial y^{\prime }}+\mathbf{i}%
_{z}^{\prime }\frac{\partial H}{\partial z^{\prime }}\mathbf{,\quad }K=-%
\frac{\partial H}{\partial y^{\prime }},\quad \mathbf{K}^{\prime }=\mathbf{i}%
_{x}^{\prime }\frac{\partial ^{2}H}{\partial x\partial y^{\prime }}-\mathbf{i%
}_{y}^{\prime }\frac{\partial ^{2}H}{\partial y^{\prime 2}}-\mathbf{i}%
_{z}^{\prime }\frac{\partial ^{2}H}{\partial z^{\prime }\partial y^{\prime }}%
.  \label{4.2}
\end{equation}%
The above expressions are obtained for $h\geqslant 0$. For negative $h$ the
following symmetry relations can be applied 
\begin{eqnarray}
I_{j}\left( x^{\prime },-y^{\prime },z^{\prime }\right)  &=&H\left(
l_{j}-x^{\prime },y^{\prime },z^{\prime }\right) -H\left( -x^{\prime
},y^{\prime },z^{\prime }\right) ,  \label{4.3} \\
\mathbf{I}_{j}^{\prime }\left( x^{\prime },-y^{\prime },z^{\prime }\right) 
&=&\widehat{\mathbf{H}}^{\prime }\left( l_{j}-x^{\prime },y^{\prime
},z^{\prime }\right) -\widehat{\mathbf{H}}^{\prime }\left( -x^{\prime
},y^{\prime },z^{\prime }\right) ,\quad   \nonumber \\
J_{j}\left( x^{\prime },-y^{\prime },z^{\prime }\right)  &=&\widehat{K}%
\left( l_{j}-x^{\prime },y^{\prime },z^{\prime }\right) -\widehat{K}\left(
-x^{\prime },y^{\prime },z^{\prime }\right) ,\quad   \nonumber \\
\mathbf{J}_{j}^{\prime }\left( x^{\prime },-y^{\prime },z^{\prime }\right) 
&=&\widehat{\mathbf{K}}^{\prime }\left( l_{j}-x^{\prime },y^{\prime
},z^{\prime }\right) -\widehat{\mathbf{K}}^{\prime }\left( -x^{\prime
},y^{\prime },z^{\prime }\right) ,\quad   \nonumber
\end{eqnarray}%
where 
\begin{equation}
\widehat{\mathbf{H}}^{\prime }=-\mathbf{i}_{x}^{\prime }\frac{\partial H}{%
\partial x}-\mathbf{i}_{y}^{\prime }\frac{\partial H}{\partial y^{\prime }}+%
\mathbf{i}_{z}^{\prime }\frac{\partial H}{\partial z^{\prime }}\mathbf{%
,\quad }\widehat{K}=\frac{\partial H}{\partial y^{\prime }},\quad \widehat{%
\mathbf{K}}^{\prime }=-\mathbf{i}_{x}^{\prime }\frac{\partial ^{2}H}{%
\partial x\partial y^{\prime }}-\mathbf{i}_{y}^{\prime }\frac{\partial ^{2}H%
}{\partial y^{\prime 2}}+\mathbf{i}_{z}^{\prime }\frac{\partial ^{2}H}{%
\partial z^{\prime }\partial y^{\prime }}.  \label{4.4}
\end{equation}

\subsubsection{Primitives}

The use of the BEM with center panel collocation presumes that the triangle
size is much smaller than the wavelength, so $k\left| \mathbf{r}-\mathbf{r}%
_{j}\right| \ll \pi .$ The triangle inequality then shows 
\begin{equation}
k\left| r-r_{0}\right| =k\left| \left| \mathbf{y}-\mathbf{x}\right| -\left| 
\mathbf{y}-\mathbf{x}_{j}\right| \right| \leqslant k\left| \mathbf{x}-%
\mathbf{x}_{j}\right| \ll \pi ,\quad r_{0}=\sqrt{x^{\prime 2}+y^{\prime
2}+z^{\prime 2}}.  \label{5}
\end{equation}%
For any integration point, we can use the following truncated expansion of
the exponent 
\begin{equation}
e^{ikr}=e^{ikr_{0}}e^{ik\left( r-r_{0}\right) }=e^{ikr_{0}}\sum_{m=0}^{p-1}%
\frac{\left( ik\right) ^{m}}{m!}\left( r-r_{0}\right) ^{m}+O\left( \frac{1}{%
p!}\max \left( k^{p}\left| r-r_{0}\right| ^{p}\right) \right) .  \label{6}
\end{equation}%
Using the binomial formula for $\left( r-r_{0}\right) ^{m}$, substituting
expansion into Eq. (\ref{4}), and neglecting the truncation error term, we
obtain 
\begin{equation}
H\left( x,y^{\prime },z^{\prime }\right) =\frac{1}{4\pi ik}\left[
e^{iky^{\prime }}k_{0}\left( x,y^{\prime },z^{\prime }\right)
-e^{ikr_{0}}\sum_{l=0}^{p-1}A_{l}^{(p)}k_{l}\left( x,y^{\prime },z^{\prime
}\right) \right] ,  \label{H4.3}
\end{equation}%
\begin{equation}
A_{l}^{(p)}=\frac{\left( ik\right) ^{l}}{l!}a_{p-l}\left( -ikr_{0}\right)
,\quad a_{l}\left( \xi \right) =\sum_{m=0}^{l-1}\frac{\xi ^{m}}{m!}.\quad
l=0,...,p-1,  \label{H4.4}
\end{equation}%
and the functions 
\begin{equation}
k_{m}\left( x,y^{\prime },z^{\prime }\right) =z^{\prime }\int \frac{r^{m}}{%
x^{2}+z^{\prime 2}}dx,\quad r=\sqrt{x^{2}+y^{\prime 2}+z^{\prime 2}},\quad
m=0,\pm 1,\pm 2,...,  \label{H4.5}
\end{equation}%
can be recursively computed as described in the Appendix.

Computation of the derivatives of the single and double layer potentials
requires derivatives of the primitives with respect to any of their three
arguments. This can be done using a similar method, with the only notice,
that, first, the derivative should be computed using integral form (\ref{4})
and then truncated expansion (\ref{6}) should be inserted there. Dropping
these derivations we present only the final results, which can be written as
follows (for all functions arguments are $\left( x,y^{\prime },z^{\prime
}\right) $). 
\begin{eqnarray}
H &=&\frac{1}{ik}\left( R-S\right) ,\quad \frac{\partial H}{\partial x}%
=-z^{\prime }T,\quad \frac{\partial H}{\partial y^{\prime }}=R-V,\quad \frac{%
\partial H}{\partial z^{\prime }}=xT-U,  \label{H4.6} \\
\frac{\partial ^{2}H}{\partial x\partial y^{\prime }} &=&-z^{\prime }P,\quad 
\frac{\partial ^{2}H}{\partial y^{\prime 2}}=-k^{2}H+z^{\prime }Q,\quad 
\frac{\partial ^{2}H}{\partial z^{\prime }\partial y^{\prime }}=xP-y^{\prime
}Q,  \nonumber
\end{eqnarray}%
\begin{eqnarray}
P &=&\frac{1}{4\pi }\frac{1}{x^{2}+z^{\prime 2}}\left( \frac{y^{\prime }}{r}%
e^{ikr}-e^{iky^{\prime }}\right) ,\quad R=\frac{1}{4\pi }e^{iky^{\prime
}}k_{0},\quad T=\frac{1}{4\pi ik}\frac{e^{ikr}-e^{iky^{\prime }}}{%
x^{2}+z^{\prime 2}},  \label{H4.7} \\
S &=&\frac{1}{4\pi }e^{ikr_{0}}\sum_{l=0}^{p-1}A_{l}^{(p)}k_{l},\quad U=%
\frac{1}{4\pi }e^{ikr_{0}}\sum_{l=0}^{p-1}A_{l}^{(p)}i_{l-1},\quad V=\frac{%
y^{\prime }}{4\pi }e^{ikr_{0}}\sum_{l=0}^{p-1}A_{l}^{(p)}k_{l-1},  \nonumber
\\
Q &=&\frac{1}{4\pi }e^{ikr_{0}}\sum_{l=0}^{p-1}A_{l}^{(p)}\left(
iki_{l-2}-i_{l-3}\right) ,\quad   \nonumber
\end{eqnarray}%
where functions 
\begin{equation}
i_{m}\left( x;a\right) =\int r^{m}dx,\quad r=\sqrt{x^{2}+a^{2}},\quad
a^{2}=y^{\prime 2}+z^{\prime 2},\quad m=0,\pm 1,\pm 2,...,  \label{H4.8}
\end{equation}%
can be recursively computed as described in the Appendix.

We note now that in computation of primitives some cases with removable or
non-removable (strong) singularities can appear. Strong singularities appear
only for the cases when the evaluation point belongs to contour $C$, in
which case, indeed the hypersingular integral blows up (in fact,
cancellation of these singularities happens when summing up integrals over
neighbor elements; such situation is typical for evaluation of integrals
which exist in a sense of principal Cauchy value, but blow up if the
integration bounds coincide with singular points). However, this case is
never encountered in the center panel collocation method, as the center of
the element never belongs to its edge. Field points also are non-singular
since they do not belong to the boundary. However, such singularities should
be avoided in evaluation of spatial gradients of the double layer potential
(e.g. evaluation point should not approach the boundary closer than some
small distance consistent with the flat constant panel discretization errors.

Note that all singularities occur only if $z^{\prime }=0$. In this case,
however, 
\begin{equation}
H=\frac{\partial H}{\partial x}=\frac{\partial H}{\partial y^{\prime }}=%
\frac{\partial ^{2}H}{\partial x\partial y^{\prime }}=\frac{\partial ^{2}H}{%
\partial y^{\prime 2}}=0,  \label{H4.9}
\end{equation}%
so these functions should be set to zero immediately. Other derivatives are
regular and computable if $x\neq 0$ and $y^{\prime }\neq 0.$ Strong singular
case corresponds to $y^{\prime }=0$ and $x=0$, which should be avoided as
the program should avoid computation of integrals at edge evaluation points
where $y^{\prime }=0$ and $x=0$ appears on the integration path. Formulae
for limiting cases $x=0$ and $y\neq 0$ and $x\neq 0$ and $y=0$ can be found
straightforward from the above analytical expressions.

In case $k\rightarrow 0$ the Green function of the Helmholtz equation
smoothly transits to the Green function of the Laplace equation. In this
case the series related toe the expansions of the exponents disappear and
all primitives can be expressed via elementary functions: 
\begin{eqnarray}
H &=&\frac{y^{\prime }k_{0}+k_{1}}{4\pi },\quad \frac{\partial H}{\partial x}%
=-\frac{z^{\prime }}{4\pi }\frac{r-y^{\prime }}{x^{2}+z^{\prime 2}},\quad 
\frac{\partial H}{\partial y^{\prime }}=\frac{k_{0}-y^{\prime }k_{-1}}{4\pi }%
,  \label{H9} \\
\frac{\partial H}{\partial z^{\prime }} &=&\frac{x}{4\pi }\frac{r-y^{\prime }%
}{x^{2}+z^{\prime 2}}-\frac{i_{-1}}{4\pi },\quad \frac{\partial ^{2}H}{%
\partial x\partial y^{\prime }}=-\frac{z^{\prime }}{4\pi }\frac{1}{\left(
x^{2}+z^{\prime 2}\right) }\left( \frac{y^{\prime }}{r}-1\right) ,  \nonumber
\\
\frac{\partial ^{2}H}{\partial y^{\prime 2}} &=&-\frac{z^{\prime }i_{-3}}{%
4\pi },\quad \frac{\partial ^{2}H}{\partial z^{\prime }\partial y^{\prime }}=%
\frac{x}{4\pi }\frac{1}{\left( x^{2}+z^{\prime 2}\right) }\left( \frac{%
y^{\prime }}{r}-1\right) +\frac{y^{\prime }i_{-3}}{4\pi }.  \nonumber
\end{eqnarray}%
These expressions are also useful for understanding of the singularities in
general case $k\neq 0$, as these singularities are caused by singularities
of the Laplacian kernel can be efficiently computed as described in the
Appendix.

\appendix

\section{Elementary integral computations}

\subsection{Computation of integrals $i_m$}

In this appendix we show how integrals 
\begin{equation}
i_{m}\left( x;a\right) =\int r^{m}dx,\quad r=\sqrt{x^{2}+a^{2}},\quad
a^{2}=y^{\prime 2}+z^{\prime 2},\quad m=0,\pm 1,\pm 2,...  \label{A1}
\end{equation}
can be computed analytically. For small values of even and odd $m$ these
integrals can be computed 
\begin{eqnarray}
i_{0}\left( x;a\right) &=&\int dx=x,\quad i_{-1}\left( x;a\right) =\int
r^{-1}dx=\ln \left| r+x\right| ,  \label{A1.1} \\
i_{-2}\left( x;a\right) &=&\int r^{-2}dx=\frac{1}{a}\arctan \frac{x}{a}%
,\quad i_{-3}\left( x;a\right) =\int r^{-3}dx=\frac{x}{a^{2}r}.  \nonumber
\end{eqnarray}
We have the recurrence: 
\begin{eqnarray}
i_{m+2} &=&\int r^{m+2}dx=xr^{m+2}-\left( m+2\right) \int x^{2}r^{m}dx
\label{A1.2} \\
&=&xr^{m+2}-\left( m+2\right) i_{m+2}+(m+2)a^{2}i_{m}.  \nonumber
\end{eqnarray}
\begin{equation}
i_{m+2}=\frac{xr^{m+2}}{m+3}+\frac{m+2}{m+3}a^{2}i_{m},\quad m\neq -3.
\label{A1.3}
\end{equation}
The special case $m=-3$ is not so important, since we know for $i_{-1}\left(
x;a\right) .$ We need integral values for $m\geqslant -3$, and have explicit
expressions for all required non-positive $m $ and we can recursively find
all positive $m$ starting the recurrence from $m=0$ for even $m$ and from $%
m=-1$ for odd $m$.  For even $m=2n$ we can evaluate the integrals as series
explicitly using the binomial expansion 
\begin{equation}
\!\!\!\!\!\!\!\! i_{2n}\left( x;a\right) =\int r^{2n}dx=\int \left( x^{2}+a^{2}\right)
^{n}dx=\sum_{l=0}^{n}\frac{n!}{l!(n-l)!}a^{2(n-l)}\int
x^{2l}dx=\sum_{l=0}^{n}\frac{n!a^{2(n-l)}x^{2l+1}}{l!(n-l)!\left(
2l+1\right) }.  \label{A2}
\end{equation}
For odd $m=2n+1$ we have the recurrence 
\begin{equation}
i_{2n+1}=\frac{1}{2n+2}\left[ xr^{2n+1}+\left( 2n+1\right) a^{2}i_{2n-1}%
\right] .  \label{A5}
\end{equation}
Expanding this down to $i_{-1}$ we obtain 
\begin{equation}
i_{2n+1}\left( x;a\right) =x\sum_{l=0}^{n}\frac{2^{2l}\left( 2n+2\right)
!(l!)^{2}}{2^{2n+2}\left( 2l+1\right) !\left[ \left( n+1\right) !\right] ^{2}%
}r^{2l+1}a^{2(n-l)}+\frac{\left( 2n+2\right) !}{2^{2n+2}\left[ \left(
n+1\right) !\right] ^{2}}a^{2n+2}\ln \left| r+x\right| .  \label{A6}
\end{equation}
Note that primitives $i_{m}\left( x;a\right) $ may have singularities only
if $a=0.$ In such cases we have 
\begin{equation}
i_{m}\left( x;0\right) =\frac{x\left| x\right| ^{m}}{m+1},\quad m\neq
-1,\quad i_{-1}\left( x;0\right) =\mbox{sgn}(x)\ln \left| x\right| .  \label{A6.1}
\end{equation}%
This shows that all functions with negative index are singular at $x=0$ and
this should be accounted for in formulae involving them.

\subsection{Computation of integrals $k_m$}

In this appendix we show how integrals 
\begin{equation}
k_{m}\left( x;y^{\prime },z^{\prime }\right) =z^{\prime }\int \frac{r^{m}}{%
x^{2}+z^{\prime 2}}dx,\quad r=\sqrt{x^{2}+y^{\prime 2}+z^{\prime 2}},\quad
m=0,\pm 1,\pm 2,...  \label{B1}
\end{equation}%
can be computed analytically. Note, that for $y^{\prime }=0$ these reduce to
the integrals $i_{m}$: 
\begin{equation}
k_{m}\left( x;0,z^{\prime }\right) =z^{\prime }\int \frac{r^{m}}{r^{2}}%
dx=z^{\prime }i_{m-2}\left( x;\left| z\right| \right) .  \label{B1.0}
\end{equation}%
So we only need consider the case $y^{\prime }\neq 0.$ We derive recurrences
for all needed $k_{m}$ and find initial values to use them. We have 
\begin{eqnarray}
k_{m+2} &=&z^{\prime }\int \frac{r^{m+2}}{x^{2}+z^{\prime 2}}dx=z^{\prime
}\int \frac{\left( x^{2}+z^{\prime 2}\right) +y^{\prime 2}}{x^{2}+z^{\prime
2}}\rho ^{m}dx  \label{B1.1} \\
&=&z^{\prime }i_{m}\left( x;\sqrt{y^{\prime 2}+z^{\prime 2}}\right)
+y^{\prime 2}k_{m}.  \nonumber
\end{eqnarray}%
With known expressions for the initial values (odd and even), all integrals
can be computed recursively. For this we have 
\begin{eqnarray}
k_{0}\left( x;y^{\prime },z^{\prime }\right)  &=&\frac{z^{\prime }}{\left|
z^{\prime }\right| }\arctan \frac{x}{\left| z^{\prime }\right| }=\mbox{sgn}(z^{\prime })\arctan \frac{x}{\left| z^{\prime }\right| }  \label{B1.2}
\\
k_{1}\left( x;y^{\prime },z^{\prime }\right)  &=&\frac{y^{\prime }z^{\prime }%
}{\left| z^{\prime }\right| }\arctan \frac{y^{\prime }x}{\left| z^{\prime
}\right| r}+z^{\prime }\ln \left| r+x\right| =y^{\prime }\mbox{sgn}(z^{\prime
})\arctan \frac{y^{\prime }x}{\left| z^{\prime }\right| r}+z^{\prime }\ln
\left| r+x\right| .  \nonumber
\end{eqnarray}%
For the boundary integrals we need only one negative value $m=-1.$ We can
find this value from the same recurrence and known $k_{1}:$ 
\begin{eqnarray}
k_{-1}\left( x;y^{\prime },z^{\prime }\right)  &=&\frac{1}{y^{\prime 2}}%
\left[ k_{1}\left( x;y^{\prime },z^{\prime }\right) -z^{\prime }i_{-1}\left(
x;\sqrt{y^{\prime 2}+z^{\prime 2}}\right) \right]   \label{B1.3} \\
&=&\frac{z^{\prime }}{y^{\prime }\left| z^{\prime }\right| }\arctan \frac{%
y^{\prime }x}{\left| z^{\prime }\right| r}=\frac{\mbox{sgn}(z^{\prime })}{y^{\prime
}}\arctan \frac{y^{\prime }x}{\left| z^{\prime }\right| r}.  \nonumber
\end{eqnarray}%
Note that $y^{\prime }k_{-1}\left( x;y^{\prime },z^{\prime }\right) $
entering the primitive expressions are not singular even when $z^{\prime }$
and $y^{\prime }$ approaches zero, while depend on the path (ratio $%
y^{\prime }/z^{\prime }$). 

\bibliography{RoomAcoustics}

\begin{thebibliography}{GAD13}

\bibitem[GAD13]{Gumerov2013:POMA}
N~Gumerov, R~Adelman, and R~Duraiswami.
\newblock Fast multipole accelerated indirect boundary elements for the
  {H}elmholtz equation.
\newblock {\em Proc. Meet. Acoust.}, 2013.

\bibitem[GD09]{Gumerov2009:JASA}
Nail~A Gumerov and Ramani Duraiswami.
\newblock A broadband fast multipole accelerated boundary element method for
  the three dimensional {H}elmholtz equation.
\newblock {\em J. Acoust. Soc. Am.}, 125(1):191--205, January 2009.

\end{thebibliography}
\bibliographystyle{alpha}
\end{document}